*Review Article*

# Approximations of Maximal and Principal Ideal

Faraj Abdunanbi

*Department Mathematics, University of Ajdabiya, Libya.*

*Cite this article:* Abdunanbi F. Approximations of Maximal and Principal Ideal. Alq J Med App Sci. 2021;4(1):7-11.

## INTRODUCTION

The rough set theory has shown by Pawlak [1] in 1982. It is a good formal tool for modeling and processing incomplete information in information system. In recently 40 years, some researchers develop this theory and use it in many areas. The upper approximation of a given set is the union of all the equivalence classes, which are subsets of the set, and the lower approximation is the union of all the equivalence classes, which are intersection with set non-empty. Many researchers develop and use the rough theory in the group and ring theory. For example, the notation of rough subring with respect ideal has presented by B. Davvaz [2]. Algebraic properties of rough sets have been studied by Bonikowaski [3]. John N. Mordeson [4], has used covers of the universal set to defined approximation operators on the power set of the given set. Some concept lattice in Rough set theory has studied by Y.Y. Yao [5]. Ronnason Chinram, [6] studies rough prime ideas and Rough fuzzy prime ideals in gamma-semigroups. Some other substitute an algebraic structure instead of the universe set. Like Biswas and Nanda [7], they make some notions of rough subgroups. Kuroki in [8], introduced the notion of a rough ideal in a semi group. Some properties of the Also, Kuroki and Mordeson in [9] studied the structure of rough sets and rough groups. S. B Hosseinin [10], he introduced and discussed the concept of T-rough (prime, primary) ideal and T-rough fuzzy (prime, primary) ideal in a commutative ring. In addition, B. Davvaz [11] applied the concept of approximation spaces in the theory of algebraic hyperstructures, and in investigated the similarity between rough membership functions and conditional probability. In this paper, we shall introduce the maximal ideal. Our result will introduce the rough maximal ideal as an extended notion of a classic maximal ideal and we study some properties of the lower and the upper approximations a maximal ideal.

## PRELIMINARIES

Suppose that $U$ (*universe*) be a nonempty finite set. Let $R$ *an equivalence relation* (reflexive, symmetric, and transitive) on an $U$. Some authors say $R$ is *indiscernibility relation*. The discernibility relation represents our lack of knowledge about elements of $U$. For simplicity, they assume $R$ *an equivalence relation*. We use $U/R$ to denote the family of all equivalent classes of $R$ (or classifications of $U$), and we use $[x]_R$ to denote an equivalence class in $R$ containing an element $x \in U$. The pair $(U, R)$ is called an approximation space. The empty set $\varnothing$ and the elements of $U/R$ are call elementary sets. For any $X \subseteq U$, we write $X^c$ to denote the complementation of $X$ in $U$.

**Definition 2.1**: For an approximation space $(U, R)$, we define the upper approximation of $X$ by $\overline{RX} = \{x \in U : [x]_R \cap X \neq \emptyset\}$; i.e. $\overline{RX}$ the set of all objects which can be only classified as *possible* members of $X$ with respect to $R$ is called the $R$-*upper approximation* of a set $X$ with respect to $R$. And the lower approximation of $X$ by $\underline{RX} = \{x \in U : [x]_R \subseteq X\}$. i.e $\underline{RX}$ is the set of all objects which can be with *certainty* classified as members of $X$ with respect to $R$ is called the $R$-*lower approximation* of a set $X$ with respect to $R$.

**Definition 2.2**: For an approximation space $(U, R)$, we define the *boundary region* by $BX_R = \overline{RX} - \underline{RX}$.
i.e. $BX_R$ is the set of all objects which can be decisively classified neither as members of $X$ nor as the members of $X^c$ with respect to $R$. If $BX_R = \emptyset$, we say $X$ is exact (*crisp*) set. But if $BX_R \neq \emptyset$, we say $X$ *Rough set* (inexact). We have Properties of approximations as:

1) $\underline{RX} \subseteq X \subseteq \overline{RX}$
2) $\underline{R\emptyset} = \overline{R\emptyset}, \underline{RU} = \overline{RU}$,
3) $\underline{R(X \cup Y)} \supseteq \underline{R(X)} \cup \underline{R(Y)}$,
4) $\underline{R(X \cap Y)} = \underline{R(X)} \cap \underline{R(Y)}$,
5) $\overline{R(X \cup Y)} = \overline{R(X)} \cup \overline{R(Y)}$
6) $\overline{R(X \cap Y)} \subseteq \overline{R(X)} \cap \overline{R(Y)}$
7) $\overline{RX^C} = (\underline{RX})^c$
8) $\underline{RX^C} = (\overline{RX})^c$
9) $\underline{R(\underline{RX})} = \overline{R(\underline{RX})} = \underline{RX}$
10) $\overline{R(\overline{RX})} = \underline{R(\overline{RX})} = \overline{RX}$

**Example 2.1:** Let us consider set of objects $U = \{x_1, x_2, x_3, x_4, x_5, x_6, x_7\}$, and the Equivalence relation $R = \{\{x_1\}, \{x_2\}, \{x_3, x_4\}, \{x_5, x_7\}, \{x_6\}\}$, and let $X = \{x_1, x_4, x_6\}$. Now, the upper approximations of X: $\overline{RX} = \{x \in U : [x]_R \cap X \neq \emptyset\}$. So $\overline{RX} = \{x_1, x_3, x_4, x_6\}$ and the lower approximation of X : $\underline{RX} = \{x \in U : [x]_R \subseteq X\}$. So $\underline{RX} = \{x_1, x_6\}$. The *boundary region* by $BX_R = \overline{RX} - \underline{RX}$. So, $BX_R = \{x_3, x_4\}$. Then $BX_R \neq \emptyset$, **so** X is *Rough set*.

Now, we define the ideal and maximal idea. Then we will study the upper and lower approximations ideal. We suppose we have a ring $\mathcal{R}$ and $I$ be an Ideal of a ring $\mathcal{R}$, and $X$ be a non-empty subset of $\mathcal{R}$.

**Definition 2.3**: Let $I$ be an Ideal of $\mathcal{R}$; For $a, b \in \mathcal{R}$ we say $a$ is congruent of $b$ mod I, we express this fact in symbols as $a \equiv b \pmod{I}$ if $a - b \in I$ ………………(1)

Not that, it easy to see the relation (1) is an equivalents relation. Therefore, when we let $U = \mathcal{R}$ and we suppose a relation $R$ is the equivalents relation (1), so we can defined the upper and lower approximation of X with respect of I as: $\overline{I(X)}$ =∪ { $x \in \mathcal{R} : (x + I) \cap X \neq \emptyset$} , $\underline{I(X)}$ =∪ { $x \in \mathcal{R} : x + I \subseteq X$} , respectively. Moreover, the boundary of X with respect of I is B$X = \overline{I(X)} - \underline{I(X)}$. If BX=∅ we say X is *Rough set with respect I.* For any approximation space ($U,R$) by rough approximation on ($U,R$), we mean a mapping Apr(X): $p(U) \to P(U) \times P(U)$ defined by for all $x \in P(U)$, Apr(X)= $(\overline{I(X)}, \underline{I(X)})$, where $\overline{I(X)}$ = { $x \in \mathcal{R} : (x + I) \cap X \neq \emptyset$}, $\underline{I(X)}$ = { $x \in \mathcal{R} : x + I \subseteq X$}.

**Definition 2.4**: An ideal $M$ in a ring $\mathcal{R}$ we called maximal if M ≠ $\mathcal{R}$ and the only ideal strictly containing $M$ is $\mathcal{R}$.

**Definitions 2.5:** Let $\mathcal{R}$ be a commutative ring with identity. Let S be a subset of $\mathcal{R}$. The ideal generated by S is the subset < S > = {$r_1s_1 + r_2s_2 + \ldots + r_ks_k \in \mathcal{R}$ | $r_1, r_2, \ldots \in \mathcal{R}$, $s_1, s_2, \ldots \in S, k \in N$}. In particular, if S has a single element s this is called the principal ideal generated by s. That is, < s > = {$rs \mid r \in \mathcal{R}$ }.

**Examples 2.2:** The ideal 2$\mathbb{Z}$ of $\mathbb{Z}$ is the principal ideal <2>.

**Examples 2.3:** In $\mathbb{Z}$, the ideal < 5 > is maximal. For suppose that $I$ is an ideal of $\mathbb{Z}$ properly containing <5>. Then there exists some $m \in I$ with $m \notin$ < 5 >, i.e. 5 does not divide $m$. Then gcd(5, $m$) = 1 since 5 is prime, and we can write 1 = 5$x$ + $my$ for integers $x$ and $y$. Since 5$x \in I$ and $my \in I$, this means 1 ∈ $I$. Then $I = \mathbb{Z}$, and < 5 >, is a maximal ideal in $\mathbb{Z}$. Note that the maximal ideals in $\mathbb{Z}$ are precisely the ideals of the form < $p$ >, where $p$ is prime.

**Proposition 1.4** Let $\mathcal{R}$ be a commutative ring with identity. Then every maximal ideal of $\mathcal{R}$ is prime.

**Example 2.4**. For $\mathcal{R} = \mathbb{Z}_{12}$, two maximal ideals are $M_1$ = {0, 2, 4, 6, 8, 10} and $M_2$ = {0, 3, 6, 9}. Two other ideals, which are not maximal are {0, 4, 8} and {0, 6}.

## UPPER AND LOWER MAXIMAL IDEALS

In this section, we introduce the maximal ideal and we study some properties of upper and lower maximal ideal. Let consider the example2-4:

**Example 3.1**. Let consider the example2-4 $\mathcal{R} = \mathbb{Z}_{12}$, $M$ = {0, 3, 6, 9}. X={1,2,6,7,9} For ∈ $\mathcal{R} : x + M$, we get {0,3,6,9}, {1,4,7,10},{2,5,8,11}. Now, the upper approximations of X with respect of M: $\overline{M(X)}$ =∪ { $x \in \mathcal{R} : (x + M) \cap X \neq \emptyset$}, = {0,1,2,3,4,5,6,7,8,9,10,11} and lower approximation of X with respect of M= $\underline{M(X)}$ = ∪ { $x \in \mathcal{R} : x + M \subseteq X$}, So, $\underline{I(X)} = \emptyset$ because no element satisfy the definition of $\underline{M(X)}$. Moreover, BX =$\overline{I(X)} - \underline{I(X)}$ = {0,1,2,3,4,5,6,7,8,9,10,11}. Then, X is rough set with respect $M$.

**Example 3.2.** Let us consider the ring $\mathcal{R} = \mathbb{Z}_6$. Suppose let maximal ideals is $M$={0,2,4} and X={1,2,3,4,5}. For ∈ $\mathcal{R} : x + M$, we get {0,2,4}, {1,3,5}. The upper approximations of X with respect of $M$: {0, 2, 4} ∪{1, 3, 5}.so, $\overline{M(X)}$ = {0,1,2,3,4,5} and the lower approximation of X with respect of $M$ :$\underline{M(X)}$ = {1,3,5}. B$X = \overline{M(X)} - \underline{M(X)}$ = {0,2,4}. Then $X$ is rough set with respect maximal ideal $M$.

We can study the properties of maximal ideal in next proposition:

**Proposition 3-1**: For every approximation ($\mathcal{R},M$) and Every subset $A, B \subseteq \mathcal{R}$ we have:
1) $\underline{M(A)} \subseteq A \subseteq \overline{M(A)}$;
2) $\underline{M(\emptyset)} = \emptyset = \overline{M(\emptyset)}$;
3) $\underline{M(\mathcal{R})} = \mathcal{R} = \overline{I(\mathcal{R})}$;
4) $\underline{M(A \cap B)} = \underline{M(A)} \cap \underline{M(B)}$;
5) $\overline{M(A \cup B)} = \overline{M(A)} \cup \overline{M(B)}$
6) If A $\subseteq$ B , then $\underline{M(A)} \subseteq \underline{M(B)}, and \overline{M(A)} \subseteq \overline{M(B)}$;
7) $\underline{M(A \cup B)} \supseteq \underline{M(A)} \cup \underline{M(B)}$;
8) $\overline{M(A \cap B)} \subseteq \overline{M(A)} \cap \overline{M(B)}$
9) $\overline{M(A)} = ( \underline{M(A^c)})^c$
10) $\underline{M(A)} = (\overline{M(A^c)})^c$
11) $\underline{M(\underline{I(A)}} = \overline{M(\underline{M(A)})} = \underline{M(A)}$
12) $( \overline{M(\overline{M(A)}} = (\overline{M(A)} = \overline{M(A)}$

*Proof:*
1) If $x \in M(A)$, then $x \in M(A) = \{x \in \mathcal{R} : x + M \subseteq A$, then $x \in A$, Hence $M(A) \subseteq A$, next if $x \in A$, $\overline{M(A)} = \{x \in \mathcal{R} : (x + M) \cap A \neq \emptyset\}$ ,, then $x \in \overline{M(A)}$ then $A \subseteq \overline{M(A)}$.
2) And 3) it easy to see that.
4) If $x \in M(A \cap B)$, then $x \in M(A \cap B) = \{x \in \mathcal{R} : x + M \subseteq A \cap B\}$, then $x \in \mathcal{R} : x + M \subseteq A$ and $x + MB$, then then $x \in M(A) \cap M(B)$.
5) It says way in 4)
6) Since $A \subseteq B$, then $A \cap B = A$, by 4) then $M(A) \cap M(B)$; It implies $M(A) \subseteq M(B)$, also, by 5) we get $\overline{M(A)} \subseteq \overline{M(B)}$.

7) Since $A \subseteq A \cup B$, $B \subseteq A \cup B$, by 6) we get $M(A) \subseteq M(A \cup B)$, and $M(B) \subseteq M(A \cup B)$ ; wich yields $M(A \cup B) \supseteq M(A) \cup M(B)$;
8) It says way in 7)
9) -13) it is easy to see that by using of definition of upper and lower approximations of $A$ with respect $M$.

## UPPER AND LOWER PRINCIPAL IDEAL:

**Example 4.1.** Let us consider the ring $\mathcal{R} = \mathbb{Z}_6$, $J=\{0,2,4\}$ and $X=\{1,2,3,4,5\}$. For $x \in \mathcal{R} : x + J$, we get $\{0,2,4\}$, $\{1,3,5\}$. The upper approximations of $X$ with respect of $J$: $\{0, 2, 4\} \cup \{1, 3, 5\}$. so, $\overline{J(X)} = \{0,1,2,3,4,5\}$. And the lower approximation of X with respect of $J$: $J(X) = \{1,3,5\}$, $BX = \overline{J(X)} - J(X) = \{0,2,4\}$. Then $X$ is rough set with respect $J$.

**Proposition 4-1**: suppose $\mathcal{R}$ commtive ring and $I$ principle ideal. For every approximation $(\mathcal{R},I)$ w and Every subset $A, B \subseteq \mathcal{R}$ we have:
1) $I(A) \subseteq A \subseteq \overline{I(A)}$;
2) $I(\emptyset) = \emptyset = \overline{I(\emptyset)}$;
3) $I(\mathcal{R}) = \mathcal{R} = \overline{I(\mathcal{R})}$;
4) $I(A \cap B) = I(A) \cap I(B)$;
5) $\overline{I(A \cup B)} = \overline{I(A)} \cup \overline{I(B)}$
6) If $A \subseteq B$, then $I(A) \subseteq I(B)$, and $\overline{I(A)} \subseteq \overline{I(B)}$;
7) $I(A \cup B) \supseteq I(A) \cup I(B)$;
8) $\overline{I(A \cap B)} \subseteq \overline{I(A)} \cap \overline{I(B)}$
9) $\overline{I(A)} = (I(A^C))^c$
10) $I(A) = (\overline{I(A^C)})^c$
11) $I(I(A) = \overline{I(I(A))} = I(A)$

12) $\overline{(I(\overline{I(A)})} = I(\overline{I(A)}) = \overline{I(A)}$
13) $I(x + I) = \overline{I(x + I)}$ for all $x \in \mathcal{R}$.

**Definition 4.1.** If $A$ and $B$ are anon-empty subset of $\mathcal{R}$, we denote $AB$ for the set of all finite sums $\{a_1 b_1 + a_2 b_2 ,…, a_n b_n : n \in \mathbb{N}, a_i \in A, b_i \in B\}$.i.e: $AB = \sum_{i=1}^{n}(a_i \cdot b_i)$, $a_i \in A, b_i \in B$.

**Proposition 4-2:** Let $I$ be maximal or principal Ideal of $\mathcal{R}$, and $A, B$ are non-empty subset of the ring $\mathcal{R}$, then
1) $\overline{I(A)} + \overline{I(B)} = \overline{I(A + B)}$
2) $\overline{I(A.B)} = \overline{I(A)} \cdot \overline{I(B)}$
3) $I(A) + I(B) \subseteq I(A + B)$
4) $I(A.B) \supseteq I(A).I(B)$

*Proof*

1) We need to proof $\overline{I(A + B)} \subseteq \overline{I(A)} + \overline{I(B)}$ and $\overline{I(A + B)} \supseteq \overline{I(A)} + \overline{I(B)}$ So, suppose $x \in \overline{I(A + B)}$, by definition of upper approximation of $A+B$ with respect $I$, $(x+I) \cap (A+B) \neq \emptyset$. Hence there exists $s \in (x+I)$ and $y \in A+B$, also, $s = \sum_{i=1}^{n} a_i + b_i$ for some $a_i \in A, b_i \in B$. we have, $x \in s + I = \sum_{i=1}^{n}(a_i + b_i) + I = \sum_{i=1}^{n}(a_i + I) + (b_i + I)$. Then there exist $x_i \in (a_i + I)$, and $s_i \in (b_i + I)$ such that $x = \sum_{i=1}^{n} x_i + s_i$ So, $x_i \in \overline{I(A)}$ And $s_i \in \overline{I(B)}$. Because, $a_i \in (x_i + I) \cap A$ and $b_i \in (s_i + I) \cap B$, Then, $\overline{I(A + B)} \subseteq \overline{I(A)} + \overline{I(B)}$.

The other side, we suppose $x \in \overline{I(A)} + \overline{I(B)}$, Then

$$x = \sum_{i=1}^{n} a_i + b_i \text{ for some } x_i \in \overline{I(A)} \text{ And } s_i \in \overline{I(B)}.$$

Hence, $(a_i+I) \cap A \neq \emptyset$ and $s_i \in (b_i+I) \cap B \neq \emptyset$ for $1 \leq i \leq n$.
Also, $\sum_{i=1}^{n} x_i + s_i \in A + B$ and $\sum_{i=1}^{n} x_i + s_i \in \sum_{i=1}^{n}(a_i + b_i) + I$. thuse $(x + I) \cap (A + B) \neq \emptyset$, that mean $x \in \overline{I(A + B)}$, so, $\overline{I(A + B)} \supseteq \overline{I(A)} + \overline{I(B)}$.

 2) Similar 1) by using definition of AB.

3) &4) Similar way in (1)&(2) by using definition of lower approximation.

**Example 4.3.** Let consider the ring $\mathcal{R} = \mathbb{Z}_6$, $I=\{0,2,4\}$ and $A=\{1,2,3,4,5\}$, $B=\{0,1,2,4\}$, then $AB=\sum_{i=1}^{n}(a_i \cdot b_i)$, $a_i \in A$, $b_i \in B$. $AB= \{0,1,2,3,4,5\}$.

We get $\overline{I(A)} = \{0,1,2,3,4,5\}$. In addition, $\overline{I(B)} = \{0,1,2,3,4,5\}$.

Then $\overline{I(A.B)} = \{0,1,2,3,4,5\}$.

However, $\overline{I(A)} \cdot \overline{I(B)} = \{0,1,2,3,4,,5\}$.

Then $\overline{I(A.B)} = \overline{I(A)} \cdot \overline{I(B)}$. $\underline{I(A)} = \{1,3,5\}$, $\underline{I(B)} = \{0,2,4\}$, $\underline{I(A.B)} = \{0,1,2,3,4,5\}$. $\underline{I(A)} \cdot \underline{I(B)} = \{0,2,4\}$.

So, $\underline{I(A.B)} \supseteq \underline{I(A)} \cdot \underline{I(B)}$.

## CONCLUSION

In this paper, we find there exists rough maximal ideal and principle that will be as an extension of the notion of a maximal ideal and principle ideal respectively

*Acknowledgments*

I would like to thank Tracy Tian Editorial Assistant of OJDM suggestion on this paper.

*Conflict of Interest*

There are no financial, personal, or professional conflicts of interest to declare.

## REFERENCE

[1] Pawlak Z. Rough sets, Int J Inf Comp Sci. 1982;11:341–356.
[2] B. Davvaz Roughness in Rings, Inform. Sci. 2004: 164 :147-163.
[3] Z. Bonikowaski. Algebraic structures of rough sets, in: W.P. Ziarko (Ed.), Rough Sets, Fuzzy Sets and Knowledge Discovery, Springer-Verlag, Berlin. 1995: 242–247.
[4] John N. Mordeson, Elsevier. Rough set theory applied to (fuzzy) ideal theory, Fuzzy Sets and Systems.2001:121:315–324.
[5] Yao, Y.Y., Perspectives of Granular Computing Proceedings of IEEE International Conference on Granular Computing. 2005:1;85-90.
[6] Ronnason Chinram., Rough prime ideas and Rough fuzzy prime ideals in Gamma-semigroups, Commun. Korean Math. Soc. 2009:24(3):341-351.
[7] Biswas, S. Nanda. Rough groups and rough subgroups, Bull. Polish Acad. Sci. Math.1994: 42 :251–254.
[8] N. Kuroki.Rough ideals in semigroups, Inform. Sci. 1997: 100 139–163.
[9] N. Kuroki, J.N. Mordeson. Structure of rough sets and rough groups, J. Fuzzy Math. 1997:5(1):183–191.
[10] S.B. Hosseini.T- Rough (Prime, Primary) ideal and t- Rough fuzzy Prime, Primary) ideal on commutative rings, int.J, contemp. Math science. 2012;7(4):337-350.
[11] B. Davvaz. Rough sets in a fundamental ring, Bull. Iranian Math. Soc.1998: 24 :49.